\documentclass{mfatshort}
\usepackage{amsmath}
\usepackage{amssymb}

\usepackage{amscd}
\usepackage{amsthm}
\usepackage{amstext}
\usepackage{amsfonts}
\usepackage{amsopn}
\usepackage{euscript}

\newcommand{\limind}{\mathop{{\rm lim\,ind}}\limits}
\newcommand{\XX}{{\mathfrak X}}
\newcommand{\N}{{\mathbb N}}

\newcommand{\R}{{\mathbb R}}
\newcommand{\CC}{{\mathbb C}}
\newcommand{\D}{{\mathcal D}}
\newcommand{\E}{{\mathfrak E}}
\newcommand{\I}{{\mathbb I}}
\newcommand{\Imm}{{\bf\rm Im}\,}

\newcommand{\Et}{E_\theta^{(\infty)}}
\newcommand{\inti}{\int_{-\infty}^{\infty}}

\newtheorem{theorem}{Theorem}
\newtheorem{propos}{Proposition}
\newtheorem{lemma}{Lemma}
\newtheorem{nasl}{Corollary}
\theoremstyle{remark}
\date{\empty}
\newtheorem{rmk}{Remark}
\newtheorem{example}{Example}

\usepackage[english]{babel}
\begin{document}

\author{S. Torba}

\title[INVERSE THEOREMS IN THE THEORY OF APPROXIMATION ...]{INVERSE THEOREMS IN THE THEORY OF APPROXIMATION OF
VECTORS IN A BANACH SPACE WITH EXPONENTIAL TYPE ENTIRE VECTORS}

\email{sergiy.torba@gmail.com}
\address{Institute of Mathematics, National Academy of Sciences of Ukraine, 3 Tereshchenkivs'ka,
Kyiv, 01601, Ukraine}

\thanks{This work was partially supported by the Ukrainian State Foundation for Fundamental Research
(project N14.1/003).}

\date{20/09/2008}

\begin{abstract}
Arbitrary operator $A$ on a Banach space $\XX$ which is the
generator of $C_0$-group with certain growth condition at infinity
is considered. The relationship between its exponential type entire
vectors and its spectral subspaces is found. Inverse theorems on
connection between the degree of smoothness of vector $x\in \XX$
with respect to operator $A$, the rate of convergence to zero of the
best approximation of $x$ by exponential type entire vectors for
operator $A$, and the $k$-module of continuity are established.
Also, a generalization of the Bernstein-type inequality is obtained.
The results allow to obtain Bernstein-type inequalities in weighted
$L_p$ spaces.
\end{abstract}

\maketitle

\section{Introduction}
Direct and inverse theorems which establish the relationship between
the
 degree of smoothness of a function with respect to a differentiation
operator and the rate of convergence to zero of its best
approximation by trigonometric polynomials are well known in the
theory of approximation of periodic functions. Bernstein's and
Jackson's inequalities are ones among such results.

N. P. Kuptsov proposed a generalized notion of the module of
continuity, expanded onto $C_0$-groups in a Banach space
\cite{Kyptsov}. Using this notion, A. P. Terekhin \cite{Terjoxin}
proved the generalized Bernstein-type inequalities for the cases of
bounded group
 and $s$-regular group. Remind that group $\{U(t)\}_{t\in\R}$
is called $s$-regular if resolvent of its generator $A$ satisfies
condition $\exists \theta\in\R:\quad
\|R_\lambda(e^{i\theta}A^s)\|\le\frac{C}{{\mathrm Im}\lambda}$.

G. V. Radzievsky studied direct and inverse theorems
\cite{Radzievsky1997, Radzievsky1998}, using notion of
$K$-functional instead of module of continuity, but it should be
noted that $K$-functional has two-sided estimates with regard to the
module of continuity at least for bounded $C_0$-groups.

In the papers \cite{MGorbShilinst_ExpA,MGorb_OperAppr} and
\cite{Gorb_Gr_Torba} authors investigated the case of a group of
unitary operators in Hilbert space and established Bernstein-type
and Jackson-type inequalities in Hilbert spaces and their rigs.
These inequalities are used to estimate the rate of convergence to
zero of the best approximation of both finite and infinite
smoothness vectors for operator $A$
 by exponential type entire vectors.

We consider the $C_0$-groups in the Banach space, generated by the
so-called \emph{non-quasianalytic operators} \cite{LubMatsaev}, i.e.
the groups satisfying
\begin{equation}\label{NeKvaziAnalit}
\int_{-\infty}^{\infty}\frac{\ln\left\Vert U(t)\right\Vert
}{1+t^{2}}dt<\infty.
\end{equation}
We recall that the belonging of group to the $C_{0}$ class means
that for every $x\in\mathfrak {X}$ vector-function $U(t)x$ is
continuous on $\mathbb{R}$ with respect to the norm of the space
$\mathfrak {X}$.

As it was shown in \cite{MGorbShilinst_ExpA}, the set of exponential
type entire vectors for the non-quasianalytic operator  $A$ is dense
in $\XX$, so the problem of approximation by exponential type entire
vectors is correct. On the other hand, it was shown in
\cite{Gorbachuk_NeobhidnistNekvazianal} that condition
(\ref{NeKvaziAnalit}) is close to the necessary one, so in the case
when (\ref{NeKvaziAnalit}) doesn't hold, the class of entire vectors
isn't necessary dense in $\XX$, and the corresponding approximation
problem loses its meaning.

In \cite{GrushkaTorba2007} the generalized Jackson-type inequalities
for approximation by entire vectors of exponential type of
non-quasianalytic operators are established. The purpose of this
work is to obtain Bernstein-type inequalities and the analogue of
inverse theorem for such approximations, and to give some
applications of these results to weighted $L_p$ spaces. In order to
do this, it is proved that the set of  exponential type  entire
vectors of type, not exceeding some $\sigma>0$, coincides with some
spectral subspace of non-quasianalytic operator (constructed in
\cite{LubMatsaev}), and the well-developed technique for spectral
subspaces is used. The last result (coincidence of the two sets of
vectors) improves the embedding, established in
\cite{MGorbShilinst_ExpA}.

\setcounter{equation}{0}
\section{Preliminaries}
Let $A$ be a closed linear operator with dense domain of definition
$\D(A)$ in Banach space $(\XX,\left\Vert \cdot\right\Vert )$ over
the field of complex numbers.

Let $C^{\infty}(A)$ denotes the set of all infinitely differentiable
vectors of operator $A$, i.e.
\begin{equation*}
C^{\infty}(A)=\bigcap_{n\in\N_{0}}\D(A^{n}),\quad\N_{0}=\N\cup\{0\}.
\end{equation*}
For a number $\alpha>0$ we set
\begin{equation*}
\E^{\alpha}(A)=\left\{ x\in C^{\infty}(A)\,|\,\exists
c=c(x)>0\,\,\forall k\in\N_{0}\,\left\Vert A^{k}x\right\Vert \leq
c\alpha^{k}\right\}.
\end{equation*}
The set $\E^{\alpha}(A)$ is a Banach space with respect to the norm
\begin{equation*}
\left\Vert x\right\Vert
_{\E^{\alpha}(A)}=\sup_{n\in\N_{0}}\frac{\left\Vert
A^{n}x\right\Vert }{\alpha^{n}}\,.
\end{equation*}
Then $\E(A)=\bigcup_{\alpha>0}\E^{\alpha}(A)$ is a linear locally
convex space with respect to the topology of inductive limit of the
Banach spaces $\E^{\alpha}(A)$:
\begin{equation*}
\E(A)=\limind_{\alpha\rightarrow\infty}\E^{\alpha}(A).
\end{equation*}
Elements of the space $\E(A)$ are called \cite{Radyno} exponential
type entire vectors of the operator $A$. The type $\sigma(x,A)$ of
vector $x\in\E(A)$ is defined as the number
\begin{equation*}
\sigma(x,A)=\inf\left\{ \alpha>0\,:\, x\in\E^{\alpha}(A)\right\}
=\limsup_{n\rightarrow\infty}\left\Vert A^{n}x\right\Vert
^{\frac{1}{n}}.
\end{equation*}
Denote by $\Xi^\alpha(A)$ the following set
\begin{equation}\label{XiAlphaDefn}
\Xi^\alpha(A)=\big\{x\in\E(A)\,|\,\sigma(x)\le\alpha\big\}.
\end{equation}
It is easy to see that
\begin{equation}\label{XiAlphaEmbed}
    \E^\alpha(A)\subset
    \Xi^\alpha(A)=\bigcap_{\epsilon>0}\E^{\alpha+\epsilon}(A).
\end{equation}

\begin{example}\label{AinLpDef}
Let $\XX$ is one of $L_{p}(2\pi)$ ($1\leq p<\infty$) -- spaces of
integrable in $p$-th degree over $[0,2\pi]$, $2\pi$-periodical
functions or the space $C(2\pi)$ of continuous $2\pi$-periodical
functions (the norm in $\XX$ is defined in a standard way), and let
$A$ is the differentiation operator in the space $\XX$ ($\D(A)=\{
x\in\XX\cap AC(\R)\,:\, x'\in\XX\}$; $(Ax)(t)=\frac{dx}{dt}$,  where
$AC(\R)$ denotes the space of absolutely continuous functions over
$\R$). It can be proved that in such case the space $\E(A)$
coincides with the space of all trigonometric polynomials, and for
$y\in\E(A)$~ $\sigma(y,A)=\deg(y)$, where $\deg(y)$ is the degree of
the trigonometric polynomial $y$.
\end{example}

Note that all previous definitions do not change if we replace the
operator $A$ by any operator of the form $e^{i\vartheta}A,\
\vartheta\in\R$. Moreover, main results of this article -- the
theorems \ref{ThmBernsteinIneq}, \ref{InvThm} do not depend on which
operator generates the group $U(t)$ --- either $A$ or $iA$. So, in
what follows, we always assume that the operator $iA$ is the
generator of group of linear continuous operators $\{ U(t)\,:\,
t\in\mathbb{R}\}$ of class $C_{0}$ \cite{Gor1991} on $\mathfrak
{X}$. Moreover, we suppose that the operator $A$ is
non-quasianalytic.

For $t\in\mathbb{R}_{+}$, we set
\begin{equation}\label{MUDef}
M_{U}(t):=\sup_{\tau\in\mathbb{R},\,|\tau|\leq t}\left\Vert
U(\tau)\right\Vert.
\end{equation}
The estimation $\|U(t)\|\le Me^{\omega t}$ for some $M,\omega\in\R$
implies $M_{U}(t)<\infty$ (for all $t\in\R_{+}$). It is easy to see
that the function $M_{U}(\cdot)$ has the following properties:
\begin{itemize}
\item[1)] $M_{U}(t)\geq1,$ $t\in\R_{+}$;
\item[2)] $M_{U}(\cdot)$ is monotonically non-decreasing on $\R_{+}$;
\item[3)] $M_{U}(t_{1}+t_{2})\leq M_{U}(t_{1})M_{U}(t_{2})$,
$t_{1},t_{2}\in\R_{+}$.
\end{itemize}

According to \cite{Kyptsov}, for $x\in\XX$, $t\in\R_{+}$ and
$k\in\N$ we set as a generalization of module of smoothness,
\begin{align}
\label{w_{d}ef} &\omega_{k}(t,x,A)=\sup_{0\leq\tau\leq t}\left\Vert
\Delta_{\tau}^{k}x\right\Vert,\qquad\text{where}\\
 \label{Delta_{d}ef}
&\Delta_{h}^{k}=(U(h)-\I)^{k}=\sum_{j=0}^{k}(-1)^{k-j}{{j \choose
k}}U(jh),\quad k\in\N_{0},\, h\in\R\quad (\Delta_{h}^{0}\equiv 1).
\end{align}


For arbitrary $x\in\XX$ we set, according to
\cite{MGorb_OperAppr,Gorb_Gr_Torba},
\begin{equation*}
\mathcal{E}_{r}(x,A)=\inf_{y\in\Xi^r(A)}\left\Vert x-y\right\Vert
,\quad r>0,
\end{equation*}
i.e. $\mathcal{E}_{r}(x,A)$ is the best approximation of element $x$
by exponential type entire vectors $y$ of operator $A$ for which
$\sigma(y,A)\leq r$. For fixed $x$~ $\mathcal{E}_{r}(x,A)$ does not
increase and $\mathcal{E}_{r}(x,A)\rightarrow 0,\
r\rightarrow\infty$ for every $x\in\XX$ if and only if the set
$\E(A)$ of exponential type entire vectors is dense in $\XX$.
Particularly, as indicated above, the set $\E(A)$ is dense in $\XX$
if operator $A$ generates the $C_0$-group $\{ U(t)\,:\,
t\in\mathbb{R}\}$ and this group belongs to non-quasianalytic class
(that is, it satisfies \eqref{NeKvaziAnalit}).

\setcounter{equation}{0}
\section{Spectral subspaces of non-quasianalytic operators}
\label{SpectralSubspaceSection} The main instrument for proving
generalized Bernstein inequality is the theory of spectral subspaces
of non-quasianalytic operator $A$, constructed in \cite{LubMatsaev}.
Recall that spectral subspaces (denoted by $\mathcal{L}(\Delta)$)
are defined for all segments $\Delta\subset\mathbb{R}$ and are
characterized by the following properties \cite[p.446]{LubMatsaev}:
\begin{itemize}
\item[1)] The operator $A$ is defined on whole $\mathcal{L}(\Delta)$ and is bounded on it;
\item[2)] $\mathcal{L}(\Delta)$ is invariant with respect to $A$;
\item[3)] the spectrum of part $A_\Delta$ of operator $A$,
induced in $\mathcal{L}(\Delta)$, consists of intersection of
spectrum of  $A$ with the interior of segment  $\Delta$ and,
perhaps, the endpoints of segment $\Delta$. And at that, if the
endpoint of segment $\Delta$ does not belong to the spectrum of $A$,
it does not belong to the spectrum of $A_\Delta$ either;
\item[4)] if there is some subspace $\mathcal{L}$ on which the operator
$A$ is defined everywhere and is bounded, and this subspace is
invariant with respect to $A$, and at the same time the spectrum of
the $\mathcal{L}$-induced part of $A$ is included in $\Delta$, then
$\mathcal{L}\subset\mathcal{L}(\Delta)$.
\end{itemize}

Now we describe the construction of spectral subspaces and their
main properties, and later prove the relationship with the entire
vectors of exponential type. Let $\theta(t)\ (-\infty<t<\infty)$ is
the entire function of order 1 with zeroes on the positive imaginary
ray:
\begin{equation}\label{InvThetaDef}
    \theta(t)=C\prod_{k=1}^\infty\left(1-\frac
    t{it_k}\right),\quad\text{where}\ 0<t_1\le
    t_2\le\ldots,\quad\sum_{k=1}^\infty\frac 1{t_k}<\infty,
\end{equation}
$C$ is a constant. Note that $|\theta(t)|$ satisfies the conditions
$|\theta(t_{1}+t_{2})|\leq|\theta(t_{1})|\cdot|\theta(t_{2})|,\
t_{1},t_{2}\in\R$ and
$\int_{-\infty}^{\infty}\frac{\left|\ln(\alpha(t))\right|}{1+t^{2}}dt<\infty$,
{\frenchspacing i.e. }it belongs to $\mathfrak{Q}$ (for definition
of class $\mathfrak{Q}$ see \cite{GrushkaTorba2007}).

Define by $\Et$ the class of entire functions $\phi(t)$ of finite
type and order 1 which satisfies for all $m=0,1,\ldots$ and for all
$a>0$ the condition
\begin{equation}\label{InvMthetaDef}
    M_\theta^{(m,a)}(\phi):=\int_{-\infty}^\infty|t^m\theta(at)\phi(t)|\,dt<\infty.
\end{equation}
As shown in \cite[Lemma 1.1.1]{LubMatsaev}, the Fourier transform of
the functions from $\Et$ is non-quasianalytic, that is the following
property takes place:
\begin{propos}\label{InvLemma1.1.1}
For any segment $\Delta$ of real axis and for any open finite
interval $I\supset\Delta$ there exists $\phi(t)\in\Et$ such that its
Fourier transform equals one in $\Delta$ and equals zero outside
$I$.
\end{propos}
Moreover, the class $\Et$ is linear and is closed under convolutions
and differentiation.

Next step is the construction of finite functions of operator $A$.
For the $C_0$-group with non-quasianalytic generator there exists
\cite{Marchenko2} such entire function $\theta(t)$ of order 1 with
zeroes on the positive imaginary ray that
\begin{equation*}
    \|U(t)\|\le |\theta(t)|\qquad \forall t\in\R.
\end{equation*}
Lets consider arbitrary $\phi(t)\in\Et$ and construct linear
operator
\begin{equation}\label{InvPphi}
    P_\phi=\int_{-\infty}^{\infty}\phi(t)U(t)\,dt.
\end{equation}
The operator, defined by (\ref{InvPphi}), is bounded due to
(\ref{InvMthetaDef}). Next, consider arbitrary segment $\Delta$ of
the real axis and denote by $\Et(\Delta)$ the set of such functions
$\phi(t)\in\Et$ that the Fourier transform $\tilde\phi(\lambda)=1$
in some interval containing $\Delta$. Denote by
$\mathcal{L}(\Delta)$ the subspace of vectors  $x$ such that
\begin{equation}\label{InvPx=x}
    P_\phi x=x
\end{equation}
for all $\phi(t)\in\Et(\Delta)$.

Operators $P_\phi$ are useful for studying vectors $A^nx$ and for
proving of Bernstein-type inequality because of the properties
\eqref{InvPphi}, \eqref{InvPx=x} and the property
\cite[p.445]{LubMatsaev}
\begin{equation}\label{InvAPphi}
    AP_{\phi}=\overline{P_\phi A}=P_{-i\phi'},
\end{equation}
which allows to deal with derivatives of some entire functions
instead of Banach-space operators and vectors.

The following theorem shows the close relationship between spectral
subspaces and the entire vectors of exponential type.
\begin{theorem}\label{InvLemmaEmbedding}
For all $\alpha>0$
\[
\E^{\alpha}(A)\subset\Xi^\alpha(A)=\mathcal{L}([-\alpha,\alpha]),
\]
moreover, $\Xi^\alpha(A)$ is the closed subspace of $\XX$.
\end{theorem}
\begin{proof}
First we will prove the embedding
$\Xi^\alpha(A)\subset\mathcal{L}([-\alpha,\alpha])$. To do this, the
forth property of spectral subspaces (mentioned at the beginning of
this section) will be used.

Obviously, $\E^\alpha(A)$ is an invariant subspace of $A$, and so is
$\Xi^\alpha(A)$. Denote the $\Xi^\alpha$-part of $A$ as $A_\alpha$:
\begin{equation*}
    A_\alpha=A\upharpoonright\Xi^\alpha(A).
\end{equation*}
By the mentioned property of spectral subspaces, to finish the
proof, it is enough to show that
$\sigma(A_\alpha)\subset[-\alpha,\alpha]$ and that $A$ is bounded on
$\Xi^\alpha(A)$.

Lets show $\sigma(A_\alpha)\subset[-\alpha,\alpha]$. For that we
check that all points from $\CC\backslash [-\alpha,\alpha]$ are
regular.

Let $\lambda\in\R\backslash[-\alpha,\alpha]$. $\lambda$ cannot be an
eigenvalue, otherwise for some $x\in\Xi^\alpha(A)$ and for all
$n\in\mathbb{N}$ $\|A_\alpha^nx\|=|\lambda|^n\|x\|$, which implies
$x\not\in \E^{\alpha+\epsilon}(A)$ for some $\epsilon>0$, a
contradiction with (\ref{XiAlphaEmbed}). That is, $\lambda$ is not
an eigenvalue of $A_\alpha$.

The equation
\begin{equation}\label{Inv02}
    Ax-\lambda x=y
\end{equation}
has a solution
\begin{equation*}\label{Inv03}
    x=-\sum_{n=0}^\infty\frac{A^ny}{\lambda^{n+1}}.
\end{equation*}
for any $\lambda\in\R\backslash[-\alpha,\alpha]$ and
$y\in\Xi^\alpha(A)$, and this solution belongs to
 $\Xi^\alpha(A)$, so such $\lambda\in\rho(A_\alpha)$.

Let $\Imm\lambda\ne 0$. Then, as shown in \cite[p.442]{LubMatsaev},
$\lambda$ is not an eigenvalue of $A$ (as well as $A_\alpha$) and
the resolvent $R_\lambda(A)$ is defined. We set for all
$y\in\Xi^\alpha(A)\quad x=R_\lambda y$. Then
\begin{equation*}
    \|A^nx\|=\|A^nR_\lambda y\|=\|R_\lambda A^ny\|\le
    \|R_\lambda\|\cdot\|A^ny\|,
\end{equation*}
hence $x\in\Xi^\alpha(A)$ and (by definition of resolvent) $x$, $y$
satisfy the equation (\ref{Inv02}). So again
$\lambda\in\rho(A_\alpha)$.

Thus it is shown that $\{\lambda\in\R\,|\,
|\lambda|>\alpha\}\subset\rho(A_\alpha)$ and
$\{\lambda\in\CC\,|\,\Imm\lambda\ne 0\}\subset\rho(A_\alpha)$
therefore $\sigma(A_\alpha)\subset[-\alpha,\alpha]$.

To prove the boundedness of  $A$ on $\Xi^\alpha(A)$ consider the
notion of $S$-operators \cite[p.452]{LubMatsaev}\footnote{The detail
definition of $S$-operators and construction of spectral subspaces
for them goes beyond the scope of this article, thus not sited. Only
required properties are mentioned.}, it results from the following
facts (see \cite[pp.462-465 and Theorem 6.1]{LubMatsaev}):
\begin{itemize}
\item If the operator $A$ is non-quasianalytic, then it is an $S$-operator.
\item $K_\Delta^-$ is a spectral subspace $\mathcal{L}(\Delta)$ of
operator $A$.
\item There exists such bounded linear operator $\Phi_\Delta^-(A)$,
defined on the whole $\XX$, that $K_\Delta^- =
\mathrm{Ker}\,\Phi_\Delta^-(A)$.
\item Operator $A$ is defined and is bounded on whole $K_\Delta^-$.
\item If $\mathcal{L}$ is an invariant subspace of $A$ and if the spectrum of $\mathcal{L}$-induced part
$A_{\mathcal{L}}$ of operator $A$ is included into segment $\Delta$,
then $\mathcal{L}\subset K_\Delta^-$.
\end{itemize}
Moreover, from these facts it follows that $\mathcal{L}(\Delta)$ is
closed subspace. This means that closedness of $\Xi^\alpha(A)$ would
result from the first statement of theorem
($\Xi^\alpha(A)=\mathcal{L}([-\alpha,\alpha])$).

Lets prove the embedding $\mathcal{L}([-\alpha,\alpha])\subset
\bigcap_{\epsilon>0}\E^{\alpha+\epsilon}(A)=\Xi^\alpha(A)$.\footnote{This
embedding improves the result of \cite{MGorbShilinst_ExpA}: $\forall
\alpha>0\ \exists r(\alpha):\quad
\mathcal{L}([-\alpha,\alpha])\subset\E^{r(\alpha)}(A)$. }

According to
 \cite[Lemma 3.1]{GrushkaTorba2007}, there exists such entire function $K_\theta(t)$ for
 $|\theta(t)|$ that for all $r>0$ exists a constant $c_r=c_r(\theta)>0$ such that for all $z\in\CC$
\begin{equation}\label{Inv04}
    |K_\theta(rz)|\le c_r\frac{e^{r|\Imm z|}}{|\theta(|z|)|}.
\end{equation}

Lets consider $\Delta=[-\alpha,\alpha]$ and $I=(-\alpha-4\epsilon,
\alpha+4\epsilon)\supset\Delta$. According with the proof of
\cite[Lemma 1.1.1]{LubMatsaev}, the Fourier transform of a function
\begin{multline}\label{Inv05}
    \phi(t)=\frac{K_\theta^2(-\epsilon t)e^{-(\alpha+2\epsilon)it}-K^2_\theta(\epsilon
    t)e^{(\alpha+2\epsilon)it}}{-2\pi it}=\\
    =K_\theta^2(\epsilon t)\frac{e^{-(\alpha+2\epsilon)it}-e^{(\alpha+2\epsilon)it}}{-2\pi
    it}=\frac{\alpha+2\epsilon}{\pi}K_\theta^2(\epsilon
    t)\frac{\sin\big((\alpha+2\epsilon)t\big)}{(\alpha+2\epsilon)t}
\end{multline}
equals one in $\Delta$ and equals zero outside $I$. Denote by
\begin{equation*}
    \phi_{r,\epsilon}(z):=K_\theta^2(\epsilon z)\frac{\sin
    rz}{rz},\quad z\in\CC,\ r>0,\ \epsilon>0
\end{equation*}
and estimate the derivatives $\phi_{r,\epsilon}^{(n)}(t),\ t\in\R$.
Using inequality
\begin{equation*}
    \left|\frac{\sin z}{z}\right|\le\frac{\min(1,|z|)}{|z|}e^{|\Imm
    z|}\le e^{|\Imm z|}
\end{equation*}
and (\ref{Inv04}), one can found
\begin{equation}\label{Inv06}
    |\phi_{r,\epsilon}(z)|\le \frac{c_\epsilon^2 e^{2\epsilon|\Imm
    z|}}{|\theta^2(|z|)|}\cdot e^{r|\Imm z|}=\frac{c_\epsilon^2e^{(r+2\epsilon)|\Imm
    z|}}{|\theta^2(|z|)|}.
\end{equation}

Similarly to the proof of
\cite[Lemma 3.2]{GrushkaTorba2007}, Cauchy integral formula for
$\gamma_{n,r}(t):=\left\{
\zeta\in\CC\,:\,|\zeta-t|=\frac{n}{r+2\epsilon}\right\}$ and
inequality (\ref{Inv06})  allow to obtain for $t\in\R$ and $n\in\N$
\begin{multline*}
 |\phi_{r,\epsilon}^{(n)}(t)|\leq\frac{n!}{2\pi}\,\oint_{\gamma_{n,r}(t)}\frac{|\phi_{r,\epsilon}(\xi)|}{|\xi-t|^{n+1}}|d\xi|=
 \frac{n!}{2\pi}\frac{(r+2\epsilon)^{n+1}}{n^{n+1}}\,\oint_{\gamma_{n,r}(t)}|\phi_{r,\epsilon}(\xi)||d\xi|\leq\\
 \leq\frac{c^{(!)}c_\epsilon^2e^{-n}(r+2\epsilon)^{n+1}}{\sqrt{2\pi n}}
 \oint_{\gamma_{n,r}(t)}\frac{e^{(r+2\epsilon)|\Imm
    \xi-t|}}{|\theta^2(|\xi|)|}|d\xi|,
 \quad \text{where}\ c^{(!)}=\sup_{k\in\N}\,\frac{k!}{\sqrt{2\pi
 k}}\left(\frac ek\right)^k< e^{1/12}.
\end{multline*}
Using $|\theta(t+s)|\le|\theta(t)|\cdot|\theta(s)|$, it follows from
the last inequality
\begin{multline*}
|\phi_{r,\epsilon}^{(n)}(t)|\leq
\frac{c^{(!)}c_\epsilon^2e^{-n}(r+2\epsilon)^{n+1}}{\sqrt{2\pi
n}|\theta^2(t)|}
  \oint_{\gamma_{n,r}(t)}\frac{e^{(r+2\epsilon)|\Imm
    \xi-t|}\big|\theta^2\big(|(t-\xi)+\xi|\big)\big|}{|\theta^2(|\xi|)|}|d\xi|\le\\
\le\frac{c^{(!)}c_\epsilon^2e^{-n}(r+2\epsilon)^{n+1}}{\sqrt{2\pi
n}|\theta^2(t)|}
  \oint_{\gamma_{n,r}(t)}e^{(r+2\epsilon)|\Imm
    \xi-t|}\big|\theta^2(|t-\xi|)\big|\,|d\xi|\leq\\
\le c^{(!)}c_\epsilon^2\sqrt{2\pi
n}(r+2\epsilon)^{n}\bigg|\frac{\theta\big(\frac
n{r+2\epsilon}\big)}{\theta(t)}\bigg|^2.
\end{multline*}
Returning to the function $\phi(t)$ one can get
\begin{equation}\label{Inv07}
    |\phi^{(n)}(t)|=\frac{\alpha+2\epsilon}{\pi}\big|\phi_{\alpha+2\epsilon,\epsilon}^{(n)}(t)\big|\le
    \frac{c^{(!)}c_\epsilon^2\sqrt{2\pi n}}{\pi}
    (\alpha+2\epsilon)(\alpha+4\epsilon)^n \bigg|\frac{\theta\big(\frac
n{\alpha+4\epsilon}\big)}{\theta(t)}\bigg|^2.
\end{equation}

Let $x\in\mathcal{L}([-\alpha,\alpha])$. By the construction
$\phi\in\Et(\Delta)$, thus $P_\phi x=x$ and, accordingly to
(\ref{InvAPphi}),
\begin{equation*}
    \|A^nx\|=\|A^nP_\phi x\|=\|P_{(-i)^n\phi^{(n)}}x\|.
\end{equation*}
Using (\ref{InvPphi}) and (\ref{Inv07}), the following estimate for
the latter expression can be found
\begin{multline}\label{Inv08}
    \|P_{(-i)^n\phi^{(n)}}x\|\le \inti
    |\phi^{(n)}(t)\theta(t)|\,dt\cdot \|x\|\le\\
   \le\frac{c^{(!)}c_\epsilon^2\sqrt{2\pi n}}{\pi}
    (\alpha+2\epsilon)(\alpha+4\epsilon)^n \|x\| \left|\theta^2\left(\frac
n{\alpha+4\epsilon}\right)\right|\inti\frac{dt}{|\theta(t)|}.
\end{multline}
It follows from (\ref{InvThetaDef}) that
\begin{equation*}
    \inti\frac{dt}{|\theta(t)|}=c_\theta<\infty,
\end{equation*}
so there exists such $c>0$ that
\begin{equation}\label{InvBadBernstein}
    \|A^nx\|\le c\sqrt{n}(\alpha+2\epsilon)(\alpha+4\epsilon)^n \left|\theta^2\left(\frac
n{\alpha+4\epsilon}\right)\right|\|x\|,\qquad \alpha>0,\
\epsilon>0,\ n\in\N.
\end{equation}

The following relation holds
\begin{equation}\label{Inv09}
    \lim_{n\to\infty}\big(c\sqrt{n}(\alpha+2\epsilon)(\alpha+4\epsilon)^n\big)^{1/n}=\alpha+4\epsilon,\qquad
    \alpha,\epsilon\in\R_{+}.
\end{equation}
As noted in the proof of  
\cite[Theorem 3.1]{GrushkaTorba2007}, for the function $|\theta(t)|$
it holds
\begin{equation}\label{Inv10}
    \lim_{n\to\infty}\left(\left|\theta^2\left(\frac
n\alpha\right)\right|\right)^{1/n}=1,\qquad \alpha\in\R_{+},
\end{equation}
therefore from (\ref{InvBadBernstein}), (\ref{Inv09}) and
(\ref{Inv10}) one can get
\begin{equation*}
    \sigma(x,A)=\limsup_{n\rightarrow\infty}\left\Vert A^{n}x\right\Vert
^{\frac{1}{n}}\le \alpha+4\epsilon,
\end{equation*}
that is $\forall \epsilon'>0$
$x\in\E^{\alpha+4\epsilon+\epsilon'}(A)$. Due to arbitrariness of
$\epsilon$,
\begin{equation*}
    x\in\bigcap_{\epsilon>0}\E^{\alpha+\epsilon}(A),
\end{equation*}
which was to be proved.
\end{proof}

\setcounter{equation}{0}
\section{Generalized Bernstein-type inequality}
One of the well-known inequalities in approximation theory is the
Bernstein inequality. If $f(x)$ is an entire function of exponential
type $\sigma>0$, and
\begin{equation*}
|f(x)|\le M,\qquad -\infty<x<\infty,
\end{equation*}
then
\begin{equation}\label{ClassicBernstein}
    |f'(x)|\le \sigma M,\qquad -\infty<x<\infty.
\end{equation}

In this section some generalization of Bernstein inequality for
exponential type entire vectors is proved.

Note that more detail view on (\ref{InvBadBernstein}) alows to
obtain Bernstein-type inequality. Consider the relation
(\ref{Inv10}). Note that it holds uniformly for all
$\alpha\ge\alpha_0>0$. Therefore for all $\epsilon>0$ there exists
$c_\epsilon>0$ such that
\begin{equation}\label{Inv10a}
c\sqrt{n}\left|\theta^2\left(\frac
n{\alpha+4\epsilon}\right)\right|\le c_\epsilon
(1+\epsilon)^n,\qquad \forall n\in\N,\ \forall\alpha\in\R_{+}
\end{equation}
Inequalities (\ref{InvBadBernstein}) and (\ref{Inv10a}) allow to
prove

\begin{propos}
For every $\epsilon>0$ there exists $c_{\epsilon}>0$, independent of
$\alpha$ and of $n$, such that for all $\alpha>0$
\begin{equation}\label{InvBadBernstein2}
    \|A^nx\|\le
    c_{\epsilon}(1+\epsilon)^n(\alpha+2\epsilon)(\alpha+4\epsilon)^n\|x\|,\qquad
    x\in\Xi^\alpha(A)\ \text{or}\ x\in\E^\alpha(A).
\end{equation}
\end{propos}

But in contrast with the classic Bernstein inequality, the type
$\alpha$ of vector appears in (\ref{InvBadBernstein2}) in the degree
$n+1$. Lets show that the analogous inequality with the degree $n$
holds.

\begin{theorem}[Generalized Bernstein-type inequality]\label{ThmBernsteinIneq} For all vectors $x\in\E(A)$,
of type, not exceeding some $\alpha\ge 1$, the following inequality
holds
\begin{equation}\label{InvBernstein}
    \|A^nx\|\le c_n\alpha^n\|x\|,
\end{equation}
where the constants $c_n>0$ do not depend on $x$ and on $\alpha$.
\end{theorem}

\begin{proof}
Lets consider majorant $\theta(t)$ for the function $\|U(t)\|$,
constructed in \cite{Marchenko2}\footnote{The majorant is named in
\cite{Marchenko2} as $\omega(t)$, but in this article it is denoted
as $\theta(t)$ in order not to confuse it with the module of
continuity}. Remark that $\theta(t)$ is of the form
(\ref{InvThetaDef}). Similarly to the proof of theorem
\ref{InvLemmaEmbedding} and as in
 \cite[Lemma 3.1]{GrushkaTorba2007} by the function
$\theta(t)$ one can construct the entire function $K(t)$ of
exponential type.

Lets consider such function $\phi_\alpha(t)$ that its Fourier
transform equals 1 in $[-\alpha, \alpha]$ and equals 0 outside
$(-3\alpha, 3\alpha)$. According to \cite[Lemma 1.1.1]{LubMatsaev},
one can use as $\phi_\alpha(t)$ the function
\begin{equation}\label{Inv11}
    \phi_\alpha(t)=\frac{K^2\big(\frac \alpha2 t\big)\sin 2\alpha
    t}{\pi t}.
\end{equation}
Denote by
\begin{equation*}
\phi(t):=\frac{K^2\big(\frac{t}2\big)\sin 2t}{\pi t}.
\end{equation*}
Then $\phi_\alpha(t)=\alpha\phi(\alpha t)$. As it follows from
(\ref{InvPphi}) and (\ref{InvAPphi}), it is enough to estimate the
quantity
\begin{equation*}
    \inti|\phi_\alpha^{(n)}(t)\theta(t)|\,dt
\end{equation*}
to prove the theorem. For $\alpha\ge 1$ we have $|\theta(t)|\le
|\theta(\alpha t)|$ and
\begin{equation}\label{Inv12}
\inti|\phi_\alpha^{(n)}(t)\theta(t)|\,dt\le \inti|\phi^{(n)}(\alpha
t)\theta(\alpha t)|\,\alpha dt.
\end{equation}
The change of variables $\tau=\alpha\cdot t$ gives
\begin{equation*}
    \frac{d^n\phi(\alpha
    t)}{dt^n}=\frac{d^n\phi(\tau)}{d\tau^n}\cdot \alpha^n,
\end{equation*}
thus
\begin{equation*}
    \inti|\phi^{(n)}(\alpha
t)\theta(\alpha t)|\,\alpha dt = \alpha^n\cdot\inti
|\phi^{(n)}(\tau)\theta(\tau)|\,d\tau.
\end{equation*}
It is easy to see that the last integral exists and does not depend
on $\alpha$. Let it equals $c_n>0$. Then
\begin{equation*}
    \|A^n x\|=\|P_{(-i)^n\phi^{(n)}}x\|\le c_n\alpha^n\|x\|,
\end{equation*}
which was to be proved.
\end{proof}

As the consequence of theorem \ref{ThmBernsteinIneq} we get the
following estimate for an operator $\Delta_h^k$:
\begin{nasl}\label{InvLemmaBernsteinDelta}
Let $x\in\E(A)$ and $\sigma(x)\le\alpha,\ \alpha\ge 1$. Then for all
$k\in\N$
\begin{equation}\label{InvBernsteinDelta}
    \|\Delta_h^k x\|\le c_k(h\alpha)^k M_U(kh)\|x\|,
\end{equation}
where the constant $c_k$ is the same as in the theorem
\ref{ThmBernsteinIneq}, and the function $M_U(t)$ is defined by
(\ref{MUDef}).
\end{nasl}

\begin{proof}
It holds for $\Delta_h^k$:
\begin{equation*}
    \Delta_h^k x = (U(t)-\I)^kx=\int_0^t\cdots\int_0^t
    U(\xi_1+\ldots+\xi_k)A^kx\,d\xi_1\ldots d\xi_k.
\end{equation*}
By the theorem \ref{ThmBernsteinIneq},
\begin{equation*}
    \|A^kx\|\le c_k\alpha^k\|x\|,
\end{equation*}
and $\|U(\xi_1+\ldots+\xi_k)\|\le M_U(mt)$ by the definition.
Therefore,
\begin{equation*}
    \|\Delta^k_h x\|\le \int_0^t\cdots\int_0^t
    \|U(\xi_1+\ldots+\xi_k)\|\cdot \|A^kx\|\,d\xi_1\ldots d\xi_k\le
    c_k h^k M_U(kh) \alpha^k\|x\|.\qed
\end{equation*}
\renewcommand{\qed}{}
\end{proof}

\setcounter{equation}{0}
\section{Inverse theorem of approximation}
The following results generalize classical Bernstein theorem (also
known as inverse theorem).

\begin{theorem}\label{InvThm}
Let $\omega(t)$ is the function of type of module of continuity for
which the following conditions are satisfied:
\begin{enumerate}
\item $\omega(t)$ is continuous and nondecreasing for $t\in\R_{+}$.
\item $\omega(0)=0$.
\item $\exists c>0\ \forall t\in[0,1]\quad \omega(2t)\le
c\omega(t)$.
\item $\int_0^1\frac{\omega(t)}{t}dt<\infty$.
\end{enumerate}
If, for  $x\in\XX$, there exist $n\in\N$ and $m>0$ such that
\begin{equation}\label{InvEstimate}
    \mathcal{E}_r(x,A)\le \frac{m}{r^n}\omega\left(\frac
    1r\right),\qquad r\ge 1,
\end{equation}
 then $x\in\mathcal{D}(A^n)$ and for every $k\in\N$ there exists a constant
 $m_k>0$ such that
\begin{equation}\label{InvOmegaK}
    \omega_k(t,A^nx, A)\le
    m_k\left(t^k\int_t^1\frac{\omega(u)}{u^{k+1}}du+\int_0^t\frac{\omega(u)}u
    du\right),\quad 0<t\le 1/2.
\end{equation}
\end{theorem}

The following lemma is used for the proof of theorem.
\begin{lemma}\label{InvLemm}
Suppose that the function  $\omega(t)$ satisfies conditions 1~-- 3
of theorem \ref{InvThm}. If, for $x\in\XX$, there exists $m>0$ such
that
\begin{equation}\label{InvLemmEstimate}
\mathcal{E}_r(x,A)\le m\omega\left(\frac 1r\right),\qquad r\ge 1,
\end{equation}
then, for every  $k\in\N$ there exists a constant $\tilde c_k>0$
such that
\begin{equation}\label{InvLemmOmegaK}
    \omega_k(t,x,A)\le\tilde
    c_kt^k\int_k^1\frac{\omega(\tau)}{\tau^{k+1}}d\tau,\quad
    0<t\le 1/2.
\end{equation}
\end{lemma}
\begin{rmk}
As would follow from the proof, the lemma remains true under
somewhat weaker conditions than those formulated in the theorem,
namely, it is sufficient that for an element $x\in\XX$ there exist
at least one sequence $\{u_{j}\}_{j=1}^{\infty}\subset\E(A)$ such
that $\sigma(u_{j},A)\le 2^j$ and for all $j\in\mathbb N$
\begin{equation*}
\|x-u_{j}\|\le m\cdot\omega\left(\frac 1{2^j}\right).
\end{equation*}
\end{rmk}

\begin{proof}[Proof of theorem \ref{InvThm}]
As shown in the theorem \ref{InvLemmaEmbedding}, the subspaces
$\Xi^r(A)$ are closed, therefore it follows from the definition and
from (\ref{InvEstimate}) that there exists a sequence of vectors
$\{u_j\}_{j=0}^{\infty}\subset\E(A)$ such that $\sigma(u_j,A)\le
2^j$ and
\begin{equation}\label{Inv13}
    \|x-u_j\|\le\frac{m}{2^{nj}}\omega\left(\frac{1}{2^j}\right).
\end{equation}
From the inequality (\ref{Inv13}) and conditions 1, 2 one can get
$\|x-u_j\|\to 0,\ j\to\infty$, and so the vector $x$ has the
representation
\begin{equation*}
    x=u_0+\sum_{j=1}^\infty(u_j-u_{j-1}).
\end{equation*}
Due to $\sigma(u_j-u_{j-1},A)\le 2^j,\ j\in\N$, one can find from
(\ref{InvBernstein})
\begin{multline*}
    \|A^nu_j-A^nu_{j-1}\|\le c_n2^{jn}\|u_j-u_{j-1}\|\le
    c_n2^{jn}\big(\|x-u_j\|+\|x-u_{j-1}\|\big)\le \\
    \le c_n2^{jn}\left(\frac{m}{2^{nj}}\cdot\omega\left(\frac
    1{2^j}\right)+\frac {m}{2^{n(j-1)}}\cdot\omega\left(\frac
    1{2^{j-1}}\right)\right)\le
    \frac{2mc_n2^{jn}}{2^{n(j-1)}}\cdot\omega\left(\frac
    1{2^{j-1}}\right)\le\\
    \le 2^{n+1}cc_nm\cdot\omega\left(\frac 1{2^j}\right)\le
    \frac{2^{n+1}cc_nm}{\ln
    2}\int_{2^{-j}}^{2^{-j+1}}\frac{\omega(u)}u du.
\end{multline*}
Hence, $\sum_{j=1}^{\infty}(A^nu_j-A^nu_{j-1})$ is convergent. By
virtue of closedness of operator $A^n$, $x\in\mathcal{D}(A^n)$ and
\begin{equation*}
    A^n=A^nu_0+\sum_{j=1}^\infty(A^nu_j-A^nu_{j-1}),
\end{equation*}
therefore
\begin{multline*}
    \|A^nx-A^nu_{j_0}\|\le
    \sum_{j=j_0+1}^\infty\|A^nu_j-A^nu_{j-1}\|\le
    \frac{2^{n+1}cc_nm}{\ln
    2}\sum_{j=j_0+1}^\infty\int_{2-j}^{2-{j+1}}\frac{\omega(u)}u du= \\
    =\frac{2^{n+1}cc_nm}{\ln
    2}\int_0^{2^{-j_0}}\frac{\omega(u)}u du=:\tilde
    c\Omega(2^{-j_0}),\quad j_0\in\N,
\end{multline*}
where $\tilde c=\frac{2^{n+1}cc_nm}{\ln 2}$,
\begin{equation*}
    \Omega(t):=\int_0^t\frac{\omega(u)}u du.
\end{equation*}

It is easy to see that the function $\Omega(t)$ has the following
properties:
\begin{enumerate}
\item $\Omega(t)$ is continuous and monotonically nondecreasing;
\item $\Omega(0)=0$;
\item for $t\in[0,1]$, the following relation is true:
\begin{equation*}
    \Omega(2t)=\int_0^{2t}\frac{\omega(u)}u
    du=\int_0^t\frac{\omega(2u)}u du\le c\int_0^t\frac{\omega(u)}u
    du=c\Omega(t).
\end{equation*}
\end{enumerate}
Therefore, setting  $\omega(t)=\Omega(t)$ in lemma \ref{InvLemm} and
taking remark into account, we get
\begin{multline*}
    \omega_k(t,A^nx,A)\le \tilde
    c_kt^k\int_t^1\frac{\Omega(u)}{u^{k+1}}du=\frac{\tilde
    c_kt^k}{k}\left(\left.\Omega(u)\frac
    1{u^k}\right|_1^t+\int_t^1\frac{\omega(u)}{u^{k+1}}du\right)\le\\
    \le
    m_k\left(t^k\int_t^1\frac{\omega(u)}{u^{k+1}}du+\int_0^t\frac{\omega(u)}u
    du\right).
\end{multline*}
The theorem is proved.
\end{proof}

\begin{proof}[Proof of lemma \ref{InvLemm}]
By the analogy with the proof of theorem \ref{InvThm},
 it follows from (\ref{InvEstimate}) that there exists a sequence of
vectors $\{u_j\}_{j=0}^{\infty}\subset\E(A)$ such that
$\sigma(u_j,A)\le 2^j$ and
\begin{equation}\label{Inv14}
    \|x-u_j\|\le m\omega\left(\frac{1}{2^j}\right).
\end{equation}

Lets take arbitrary $h\in(0,1/2]$ and choose a number $N$ in such a
way that $\frac 1{2^{N+1}}<h\le \frac 1{2^N}$. Inequality
(\ref{Inv14}) yields
\begin{multline}\label{Inv15}
    \|u_j-u_{j-1}\|\le \|u_j-x\|+\|x-u_{j-1}\|\le\\
    \le m\omega(2^{-j})+m\omega(2^{-j+1})\le 2m\omega(2^{-j+1})\le
    2cm
    \omega(2^{-j}).
\end{multline}
By virtue of the monotonicity of $\omega(t)$
\begin{equation}\label{Inv16}
    2^k\int_{2^{-j}}^{2^{-j+1}}\frac{\omega(u)}{u^{k+1}}du\ge
    2^k\omega(2^{-j})\int_{2^{-j}}^{2^{-j+1}}\frac{1}{u^{k+1}}du=
    \frac{2^{kj}(2^k-1)}k \omega(2^{-j}) \ge 2^{kj}\omega(2^{-j}).
\end{equation}

Since $\sigma(u_j-u_{j-1},A)\le 2^j$ and $\sigma(u_0,A)\le 1$,
according to corollary \ref{InvLemmaBernsteinDelta}
\begin{align*}
    \|\Delta^k_hu_0\|&\le c_kh^kM_U(kh)\|u_0\|,\\
    \|\Delta^k_h(u_j-u_{j-1})\|&\le
    c_kh^k(2^j)^kM_U(kh)\|u_j-u_{j-1}\|,\quad j\ge 1.
\end{align*}
Relations (\ref{Inv14}) -- (\ref{Inv16}) yield
\begin{equation*}
    \|\Delta_h^k(u_j-u_{j-1})\|\le 2\tilde ch^k(2^j)^k\omega(2^{-j})\le
    2^{k+1}\tilde
    ch^k\int_{2^{-j}}^{2^{-j+1}}\frac{\omega(u)}{u^{k+1}}du,
\end{equation*}
where $\tilde c = cc_kmM_U(kh)$, and
\begin{multline*}
    \|\Delta(x-u_N)\|\le \|(U(h)-\I)^k\|\,\|x-u_N\|\le \\
    \le (M_U(h)+1)^k\|x-u_N\|\le (M_U(h)+1)^km\omega(2^{-N}).
\end{multline*}
Using these inequalities, we obtain
\begin{multline*}
    \|\Delta_h^kx\|=\bigg\|\Delta_h^ku_0+\sum_{j=1}^N\Delta_h^k(u_j-u_{j-1})+\Delta_h^k(x-u_N)\bigg\|\le
    \displaybreak[3]\\
    \le c_kM_U(kh)h^k\|u_0\|+2^{k+1}\tilde c
    h^k \sum_{j=1}^N\int_{2^{-j}}^{2^{-j+1}}
    \frac{\omega(u)}{u^{k+1}}du+(M_U(h)+1)^km\omega(2^{-N})\le
    \displaybreak[3]\\
    \le c_kM_U(kh)h^k\|u_0\|+2^{k+1}\tilde ch^k
    \int_{2^{-N}}^1\frac{\omega(u)}{u^{k+1}}du+(M_U(h)+1)^kcm\omega(h)\le
    \displaybreak[3]\\
    \le c_kM_U(kh)h^k\|u_0\|+2^{k+1}\tilde ch^k
    \int_{h}^1\frac{\omega(u)}{u^{k+1}}du+(M_U(h)+1)^kcm\omega(h)=
    \displaybreak[3]\\
    =h^k\bigg[c_kM_U(kh)\|u_0\|+2^{k+1}\tilde c\int_{h}^1\frac{\omega(u)}{u^{k+1}}du+
    (M_U(h)+1)^kcm\frac{k}{1-h^k}\int_h^1\frac{\omega(h)}{u^{k+1}}du\bigg]\le
    \displaybreak[3]\\
    \le\tilde
    c_kh^k\int_h^1\frac{\omega(u)}{u^{k+1}}du,\qquad\text{where}
\end{multline*}
\begin{equation*}
    \tilde
    c_k:=\frac{\|u_0\|c_kM_U(k/2)}{\int_{1/2}^1\frac{\omega(u)}{u^{k+1}}du}
    +2^{k+1}cc_kmM_U(k/2) +(M_U(1/2)+1)^k\frac{cmk}{1-(1/2)^k}.
\end{equation*}
The last inequality holds for all $0<h\le 1/2$. Taking into account
the definition of module of continuity (\ref{w_{d}ef}), this
inequality finishes the proof.
\end{proof}

\setcounter{equation}{0}
\section{Examples of application of abstract direct and inverse theorems in particular spaces}\label{section5}
In this section we discuss an application of the presented theory
--- the approximation of continuous functions by entire
functions in the weighted $L_p(\R, \mu^p)$ space with growing at the
infinity weight (for example, $L_1(\R,x^n)$ spaces). Similar
problems studied in several papers (see the review
\cite{Ganzburg2000}).

Lets consider the real-valued function $\mu(t)$ satisfying the
following conditions:
\begin{itemize}
\item[1)] $\mu(t)\ge 1,\quad t\in\R$;
\item[2)] $\mu(t)$ is even, monotonically non-decreasing when $t>0$;
\item[3)] $\mu(t)$ satisfies the condition
 $\mu(t+s)\le \mu(t)\cdot\mu(s),\ s,t\in\R$.
\item[4)]
$\int_{-\infty}^{\infty}\frac{\ln\mu(t)}{1+t^2}\,dt<\infty$,
\end{itemize}
or alternatively, instead of 4), the equivalent condition holds:
\begin{itemize}
\item[4')] $\sum_{k=1}^\infty\frac{\ln\mu(k)}{k^2}<\infty$.
\end{itemize}

Below are  several important classes of functions satisfying
conditions 1)--4) (see \cite{GrushkaTorba2007} for details).

1. Constant function $\mu(t)\equiv 1,\quad t\in\R$.

2. Functions with polynomial order of growth at infinity. For such
functions the following estimate holds: $\exists k\in\N,\ \exists
M\ge 1$
\begin{equation*}
\mu(t)\le M(1+|t|)^k,\quad t\in\R.
\end{equation*}

3. Functions of the form
\begin{equation*}
    \mu(t)=e^{|t|^\beta},\quad 0<\beta<1,\ t\in\R.
\end{equation*}

4. $\mu(t)$ represented as a power series for $t>0$. I.e.,
\begin{equation*}
    \mu(t)=\sum_{n=0}^\infty\frac{|t|^{n}}{m_n},
\end{equation*}
where $\{m_n\}_{n\in\N}$ is the sequence of positive real numbers
satisfying three conditions:
\begin{itemize}
\item $m_0=1$, $m_n^2\le m_{n-1}\cdot m_{n+1},\ n\in\mathbb{N}$;
\item for all $ k,l\in\mathbb{N}\quad \frac{(k+l)!}{m_{k+l}}\le
\frac{k!}{m_k}\frac{l!}{m_l}$.
\item $\sum_{n=1}^{\infty}\left(\frac
1{m_n}\right)^{1/n}<\infty$;
\end{itemize}

5. $\mu(t)$ as a module of an entire function with zeroes on the
imaginary axis. Lets consider
\begin{equation*}
    \omega(t)=C\prod_{k=1}^\infty\left(1-\frac
    t{it_k}\right),\quad t\in\R,
\end{equation*}
where $C\ge 1,\ 0<t_1\le t_2\le\ldots,\ \sum_{k=1}^\infty \frac
1{t_k}<\infty$, and set $\mu(t):=|\omega(t)|$.

Lets consider
 the space $L_p(\R,\mu^p)$ of the functions $x(s),\ s\in\R$, integrable in
$p$-th degree with the weight $\mu^p$:
\begin{equation*}
    \|x\|^p_{L_p(\R,\mu^p)}=\int_{-\infty}^{\infty}|x(s)|^p\mu^p(s)\,ds.
\end{equation*}
$L_p(\R,\mu^p)$ is the Banach space. The differential operator
\begin{equation*}
(Ax)(t)=\frac {dx}{dt},\qquad \D(A)=\{x\in L_p(\R,\mu^p)\cap
AC(\R):\ x'\in L_p(\R,\mu^p)\}.
\end{equation*}
generates the group of shifts
 $\{U(t)\}_{t\in\R}$ in the space $L_p(\R,\mu^p)$.
This group isn't bounded. As shown in \cite{GrushkaTorba2007},
\begin{equation*}
\|U(t)\|_{L_p(\R,\mu^p)}\le\mu(|t|),\quad t\in\R.
\end{equation*}

To apply the constructed theory, we need to determine how the space
$\E(A)$ and the space of exponential type entire functions are
connected. Denote by $B_\sigma$ the set of exponential functions of
entire type $\sigma$. We show that the following embedding holds
\begin{equation}\label{ExamLpEmbed1}
    \Xi^\sigma(A)\subset B_\sigma\cap L_p(\R, \mu^p).
\end{equation}

Let $f\in\Xi^\sigma(A)$. Obviously, $f\in L_p(\R, \mu^p)$. We prove
that $f\in B_\sigma$. Due to $\mu(t)\ge 1$ we have
\begin{equation*}
    \|f\|_{L_p(\R)}\le \|f\|_{L_p(\R,\mu^p)},
\end{equation*}
thus for all $n\in\N$ and for any $\epsilon >0$
\begin{equation}\label{ExamLpEmbed2}
    \|A^nf\|_{L_p(\R)}\le\|A^nf\|_{L_p(\R,\mu^p)}\le c_\epsilon(f)(\sigma+\epsilon)^n,
\end{equation}
and so we can construct a continuation of $U(t)$ onto $\CC$ by
\begin{equation*}
    U(z)=\sum_{n=0}^\infty \frac {A^n f}{n!}z^n,\quad z\in\CC.
\end{equation*}
Moreover, \eqref{ExamLpEmbed2} ensures for all $\epsilon>0$
\begin{equation*}
    \|f(x+z)\|_{L_p(\R)}=\left\|\sum_{n=0}^\infty\frac{A^n
    f}{n!}z^n\right\|\le c_\epsilon(f) \cdot
    \|f\|e^{(\sigma+\epsilon)|z|},
\end{equation*}
which means $f\in B_\sigma$, which required.

By virtue of the classical Bernstein inequality the reverse
embedding to \eqref{ExamLpEmbed1} holds for all bounded weights
$\mu(t)$. We show that it holds for all functions $\mu(t)$,
satisfying
\begin{equation*}
    \mu(t)\ge 1+R|t|
\end{equation*}
for some $R>0$ and for all $t>t_0\ge 0$. The condition on $\mu(t)$
gives us $f\in L_1(\R)$. $f\in B_\sigma$, thus it is infinitely
differentiable and by the Paley-Wiener theorem the support of its
Fourier transform is contained in $[-\sigma, \sigma]$. Lets prove
that $f\in\Xi^\sigma(A)$ by using theorem \ref{InvLemmaEmbedding}.
To do this we need to show that for all $\phi\in
E_\theta^{(\infty)}([-\sigma,\sigma])$
\begin{equation*}
    f=P_\phi f=\inti \phi(t) U(t)f dt.
\end{equation*}
Since $\phi$ is arbitrary, we can consider $\phi_1(t)=\phi(-t)\in
E_\theta^{(\infty)}([{-}\sigma,\sigma])$. Note that
\[
\inti \phi_1(t) U(t) f(x)dt = \inti \phi(t) f(x-t) dt = \phi \ast f.
\]
The Fourier transform of $\phi \ast f$ equals to
\begin{equation*}
    \widetilde{\phi\ast f} = \tilde \phi \cdot \tilde f = \tilde f,
\end{equation*}
because $\mathrm{supp}\, f\subset [-\sigma,\sigma]$, and by the
definition of $E_\theta^{(\infty)}([-\sigma,\sigma])$ we have
$\tilde\phi = 1$ on $[-\sigma, \sigma]$. Thus,
\begin{equation*}
    P_\phi f=f\qquad \forall \phi\in
    E_\theta^{(\infty)}([-\sigma,\sigma]),
\end{equation*}
so  $f\in \mathcal{L}([-\sigma,\sigma])$ and by means of theorem
\ref{InvLemmaEmbedding} $f\in \Xi^\sigma(A)$.

We have shown that $\Xi^\sigma(A)$ coincides with $B_\sigma\cap
L_p(\R, \mu^p)$. Note that $\|f-g_\sigma\|_{L_p(\R,\mu^p)}$ is
defined only for those functions that belongs to $L_p(\R,\mu^p)$
(because of $\|g_\sigma\|_{L_p(\R,\mu^p)}\le
\|f-g_\sigma\|_{L_p(\R,\mu^p)}+\|f\|_{L_p(\R,\mu^p)}$), thus the
best approximation by exponential type entire vectors is the same as
the best approximation by entire functions of exponential type.

By applying theorems \ref{ThmBernsteinIneq} and \ref{InvThm} we get
several results for the approximation theory in $L_p(\R,\mu^p)$
spaces. First two results are the direct theorems (from
\cite{GrushkaTorba2007}) for spaces $L_p(\R,\mu^p)$.

\begin{nasl}[\cite{GrushkaTorba2007}]\label{ExamLpmu2Nasl1}
For every $k\in\N$ there exists constant $\mathbf{m}_{k}(p,\mu)>0$
such that for all $ f\in L_p(\R,\mu^p)$
\begin{equation*}
\mathcal{E}_{r}(f)\leq\mathbf{m}_{k}\cdot\tilde\omega_{k}\left(\frac{1}{r},f\right),\quad
r\ge 1.
\end{equation*}
\end{nasl}

\begin{nasl}[\cite{GrushkaTorba2007}]
Let $f\in W_p^m(\R,\mu^p),\,\, m\in\N_{0}$. Then for all $k\in\N_0$
\begin{equation*}
\mathcal{E}_{r}(f)\le\mathbf{m}_{k+m}\frac{\mu\left(\frac{m}{r}\right)}{r^{m}}\widetilde{\omega}_{k}
\left(\frac{1}{r},f^{(m)}\right),\quad r\ge 1,
\end{equation*}
where constants $\mathbf{m}_{n}$ ($n\in\N$) are the same as in the
corollary \ref{ExamLpmu2Nasl1}.
\end{nasl}

\begin{nasl}
Let $f\in L_p(\R,\mu^p)\cap B_\sigma,\ \sigma\ge 1$. Then for all
$n\in\N$ there exist such constants $c_n>0$, not depending on
$\sigma$ and on $f$, that
\begin{equation*}
    \|f^{(n)}\|_{L_p(\R,\mu^p)}\le c_n\sigma^n\|f\|_{L_p(\R,\mu^p)}.
\end{equation*}
\end{nasl}

\begin{nasl}
Let $\omega(t)$ be a function of type of module of continuity for
which the following conditions are satisfied:
\begin{enumerate}
\item $\omega(t)$ is continuous and nondecreasing for $t\in\R_{+}$.
\item $\omega(0)=0$.
\item $\exists c>0\ \forall t\in[0,1]\quad \omega(2t)\le
c\omega(t)$.
\item $\int_0^1\frac{\omega(t)}{t}dt<\infty$.
\end{enumerate}
If, for $f\in L_p(\R,\mu^p)$ there exist such $n\in\N$ and $m>0$
that
\begin{equation*}
    \mathcal{E}_r(f)\le \frac{m}{r^n}\omega\left(\frac
    1r\right),\qquad r\ge 1,
\end{equation*}
then $f\in W_p^n(\R,\mu^p)$ and for every $k\in\N$ there exists such
$m_k>0$ that
\begin{equation*}
    \omega_k(t,f^{(n)})\le
    m_k\left(t^k\int_t^1\frac{\omega(u)}{u^{k+1}}du+\int_0^t\frac{\omega(u)}u
    du\right),\quad 0<t\le 1/2.
\end{equation*}
\end{nasl}

\end{document}